\documentclass[12pt]{article}   
   
\usepackage{latexsym,amsfonts,amsmath,epsfig,tabularx,amsthm,amssymb,enumerate,bm}

\usepackage{microtype}
   
\usepackage[draft=false,colorlinks,bookmarksnumbered,linkcolor=black, citecolor=black]{hyperref}

\newtheorem{theorem}{Theorem}

\newtheorem{corollary}{Corollary}

\theoremstyle{definition}
\newtheorem{example}{Example}

\renewcommand{\d}{\,\mathrm{d}}																						% Differential
\renewcommand*{\epsilon}{\varepsilon}                                   % modified epsilon
\renewcommand*{\rho}{\varrho}                                   % modified epsilon
\newcommand*{\nach}{\rightarrow}                                        % f from X \nach Y
\newcommand*{\sep}{\; \vrule \;}                                                % set seperator: { x \in X \sep x=ay }
\newcommand*{\N}{\mathbb{N}}                                            % natural numbers IN
\newcommand*{\R}{\mathbb{R}}                                            % real numbers IR
\renewcommand{\i}{\mathrm{i}}
\newcommand*{\C}{\mathbb{C}}                                            %
\newcommand*{\Z}{\mathbb{Z}}                                            %
\renewcommand*{\S}{\mathcal{S}}                                         % curved S
\newcommand*{\A}{\mathcal{A}}                                         	% curved A
\newcommand*{\algR}{\mathcal{R}}                                         	% curved R
\newcommand*{\SI}{\mathfrak{S}}                                         % fractured S
\newcommand*{\leer}{\emptyset}                                          % empty set
\newcommand*{\0}{\mathcal{O}}                                           % open set, Landau-O
\newcommand*{\id}{\mathrm{id}}                                          % identitiy 
\newcommand*{\fu}{\mathfrak{u}}																					% subset u

\newcommand*{\norm}[1]{\left\| #1 \right\|}                             % || . ||
\newcommand*{\abs}[1]{\left| #1 \right|}                                % | . |
\newcommand*{\link}[1]{(\ref{#1})}                                      % link to labeled formulas, e.g. (5)
\newcommand*{\distr}[2]{\left\langle #1, #2 \right\rangle}              % application of distributions: \distr{f}{g} -> <f,g>
\newcommand{\maxx}[1]{ \mathop{\mathrm{max}}\left\{#1\right\} }        % Min and max of sets: \max{a,b} -> max{a,b}
\newcommand{\minn}[1]{ \mathop{\mathrm{min}}\left\{#1\right\} }
\newcommand{\wrt}{w.r.t.\ }
\newcommand{\ie}{i.e.\ }
\title{On lower bounds for integration of multivariate permutation-invariant functions}     
 
\author{Markus Weimar\footnote{Philipps-University Marburg, Faculty of Mathematics and Computer Science, Hans-Meerwein-Stra{\ss}e, Lahnberge, 35032 Marburg, Germany. Email: weimar@mathematik.uni-marburg.de} \footnote{This work has been supported by Deutsche Forschungsgemeinschaft DFG (DA 360/19-1).}}
 
\begin{document}   
   
\maketitle   
   
\begin{abstract}
In this note we study multivariate integration for permutation-invariant functions from a certain
Banach space $E_{d,\alpha}$ of Korobov type in the worst case setting.
We present a lower error bound which particularly implies that in dimension~$d$ every cubature rule which reduces the initial error necessarily uses at least $d+1$ function values.
Since this holds independently of the number of permutation-invariant coordinates, 
this shows that the integration
problem can never be strongly polynomially tractable in this setting.
Our assertions generalize results due to Sloan and Wo{\'z}niakowski~\cite{SW97}.
Moreover, for large smoothness parameters $\alpha$ our bound can not be improved.
Finally, we extend our results to the case of permutation-invariant functions from Korobov-type spaces equipped with product weights.
\end{abstract} 

\section{Introduction and main result}
Consider the integration problem $\mathrm{Int}=(\mathrm{Int}_d)_{d\in\N}$,
\begin{gather*}
		\mathrm{Int}_d\colon E_{d,\alpha}\nach\C, \qquad \mathrm{Int}_d(f)=\int_{[0,1]^d} f(\bm{x})\d \bm{x},
\end{gather*}
for periodic, complex-valued functions in the Korobov class
\begin{gather*}
		E_{d,\alpha} := \left\{ f\in L_1([0,1]^d) \sep \norm{f} := \norm{f\sep E_{d,\alpha}} := \sup_{{\bm{k}}\in\Z^d} \abs{\widehat{f}({\bm{k}})} \left(\overline{k_1}\cdot\ldots\cdot\overline{k_d}\right)^\alpha < \infty \right\}
\end{gather*}
where $d\in\N$ and $\alpha>1$.
Here $\Z$ denotes the set of integers, $\N:=\{1,2,\ldots\}$, and we set $\overline{k_m}:=\maxx{1,\abs{k_m}}$.
Moreover,
for $f\in L_1([0,1]^d)$
\begin{gather*}
		\widehat{f}({\bm{k}})
		:= \distr{f}{e^{2\pi \i {\bm{k}}\, \cdot}}_{L_2} := \int_{[0,1]^d} f(\bm{x}) \, e^{-2\pi \i {\bm{k}}\bm{x}} \d\bm{x}, \qquad {\bm{k}}=(k_1,\ldots,k_d)\in\Z^d,
\end{gather*}
denotes its ${\bm{k}}$th Fourier coefficient, where $\bm{k}\bm{x}=\sum_{m=1}^d k_m\cdot x_m$, and $\i = \sqrt{-1}$.
To approximate $\mathrm{Int}_d(f)$, without loss of generality, we consider algorithms from the class of all linear cubature rules
\begin{gather}\label{quadrule}
		\A(f) := \A_{N,d}(f) := \sum_{n=1}^N w_n \, f\!\left(\bm{t}^{(n)}\right), \qquad N\in\N_0:=\{0\}\cup\N,
\end{gather}
that use at most $N$ values of the input function~$f$ at some points $\bm{t}^{(n)}\in[0,1]^d$, $n=1\ldots,N$. 
The weights $w_n$ can be arbitrary complex numbers.
Clearly, every function $f\in E_{d,\alpha}$ has a $1$-periodic extension since their
Fourier series are absolutely convergent:
\begin{gather*}
		\sum_{{\bm{k}}\in\Z^d} \abs{\widehat{f}({\bm{k}}) \, e^{2\pi \i {\bm{k}} \, \cdot}} 
		\leq \norm{f} \cdot \sum_{{\bm{k}}\in\Z^d} \left(\overline{k_1}\cdot\ldots\cdot\overline{k_d}\right)^{-\alpha}
		= \norm{f} \cdot (1+2\,\zeta(\alpha))^d < \infty.
\end{gather*}
As usual, $\zeta(s)=\sum_{m=1}^\infty m^{-s}$ is the Riemann zeta function evaluated at $s>1$. 

In \cite{SW97} Sloan and Wo{\'z}niakowski showed that for every $d\in\N$ the \emph{$N$th minimal worst case error} of $\mathrm{Int}=(\mathrm{Int}_d)_{d\in\N}$,
\begin{gather*}
			e(N,d;\mathrm{Int}_d,E_{d,\alpha}) 
			:= \inf_{\A_{N,d}} \sup_{\norm{f \sep E_{d,\alpha}}\leq 1} \abs{\mathrm{Int}_d(f) - \A_{N,d}(f)},
\end{gather*}
equals the initial error $e(0,d;\mathrm{Int}_d,E_{d,\alpha})=1$ provided that $N<2^d$.
In other words, the integration problem on the full spaces $(E_{d,\alpha})_{d\in\N}$ suffers from the 
\emph{curse of dimensionality}, since
for every fixed $\epsilon\in(0,1)$ its \emph{information complexity} 
grows exponentially with the dimension~$d$:
\begin{align*}
			n(\epsilon,d) 
			:= n(\epsilon,d;\mathrm{Int}_d,E_{d,\alpha}) 
			:= \minn{N\in\N_0 \sep e(N,d;\mathrm{Int}_d,E_{d,\alpha}) \leq \epsilon}
			\geq 2^d, \qquad d\in\N.
\end{align*}

We generalize this result to the case of 
permutation-invariant\footnote{In \cite{W12b} we used the name \emph{symmetric} what caused some confusion.} 
subspaces in the sense of~\cite{W12b}.
To this end, for $d\in\N$ let $I_d \subseteq\{1,\ldots,d\}$ be some subset of coordinates 
and consider the integration problem $\mathrm{Int}=(\mathrm{Int}_d)_{d\in\N}$ restricted to the subspace
$\SI_{I_d}(E_{d,\alpha})$ of all $I_d$-\emph{permutation-invariant} functions 
$f \in E_{d,\alpha}$.
That is, in dimension~$d$ we restrict ourselves to functions $f$ that satisfy
\begin{gather}\label{def:sym}
		f(\bm{x}) = f(\sigma(\bm{x})) \quad \text{for all} \quad \bm{x}\in[0,1]^d
\end{gather}
and any permutation $\sigma$ from
\begin{gather}\label{def:permut_set}
			\S_{I_d} := \left\{\sigma \colon \{1,\ldots,d\} \nach \{1,\ldots,d\} 
								\sep \sigma \text{ bijective and } \sigma\big|_{\{1,\ldots,d\}\setminus I_d} = \id \right\}
\end{gather}
that leaves the elements in the complement of $I_d$ fixed.
For the ease of presentation we shall use the same notation for permutations $\sigma\in\S_{I_d}$ and for the corresponding permutations $\sigma'\colon \R^d\nach\R^d$ of $d$-dimensional vectors, given by
\begin{equation*}
\bm{x} = (x_1,\ldots,x_d) \mapsto \sigma'(\bm{x}) := (x_{\sigma(1)},\ldots,x_{\sigma(d)}).
\end{equation*}
Observe that in the case $I_d=\leer$ we clearly have $\SI_{I_d}(E_{d,\alpha}) = E_{d,\alpha}$.

One motivation to study the integration problem restricted to those subspaces
is related to approximate solutions of partial differential equations.
Many approaches to obtain such solutions 
lead us to the problem of calculating high-dimensional integrals, e.g., to obtain certain wavelet coefficients.
Obviously, it is of interest whether this can be done efficiently since 
taking into account a large number of coefficients would lead to better approximations
to the exact solution.
Therefore it is important to incorporate as many structural properties 
(such as permutation-invariance of the integrands under consideration)
as possible in order to reduce the effort for every single calculation.
In information-based complexity (IBC) this effort is measured by the behavior of the information complexity
$n(\epsilon,d)$ which can be formalized by several notions of tractability.

We remind the reader that a problem is called \emph{polynomially tractable} if its information complexity 
$n(\epsilon,d)$ is bounded from above by some polynomial in $d$ and $\epsilon^{-1}$, i.e.,
\begin{gather*}
				n(\epsilon,d) \leq C\, \epsilon^{-p} \, d^q \quad \text{for some} \quad C,p>0,\, q\geq 0 \quad \text{and all} \quad \epsilon\in(0,1], \, d\in\N.
\end{gather*}
If the last inequality remains valid even for $q=0$ then we say that the problem is \emph{strongly polynomially tractable}.
Apart from that, several weaker notions of tractability were introduced recently; 
for details see, e.g., Siedlecki~\cite{P13}.
If for some fixed $s,t\in(0,1]$ the information complexity satisfies
\begin{gather}\label{def:t-weak}
		\lim_{\epsilon^{-1}+d\nach\infty} \frac{\ln n(\epsilon,d)}{\epsilon^{-s}+d^t}=0,
\end{gather}
then we have \emph{$(s,t)$-weak tractability}. 
This generalizes the well-established notion of \emph{weakly tractable} problems (which is included as the special case $s=t=1$). It is used to describe the case of at most subexponentially growing
information complexities such as, e.g., $n(\epsilon,d)=\exp(\sqrt{d})$.
Finally, a given problem is said to be \emph{uniformly weakly tractable} if the limit condition~\link{def:t-weak} holds\footnote{Note that it clearly suffices to check \link{def:t-weak} for every $s=t\in(0,1]$ in order to conclude uniform weak tractability.} for every $s,t\in(0,1]$.

Our main result states that the integration problem for permutation-invariant functions in the above sense can never be strongly polynomially tractable,
independent of the size of the sets $I_d$ which describes the number of imposed permutation-invariance conditions.
The assertion reads as follows:
\begin{theorem}\label{thm}
			Let
			\begin{gather}\label{def:N}
					N^* := N^*(d,I_d) := (\#I_d + 1) \cdot 2^{d-\# I_d},
					\qquad d\in\N.
			\end{gather}
			Then, for every $N<N^*$,
			\begin{gather}\label{bound}
					e\!\left(N,d;\mathrm{Int}_d,\SI_{I_d}(E_{d,\alpha})\right) = 1 
%					\qquad \text{for all} \quad d\in\N, \quad \alpha>1, \quad \text{and every} \quad N < N^*(d,I_d).
			\end{gather}
			and
			\begin{gather}\label{upperbound}
					e\!\left(N^*,d;\mathrm{Int}_d,\SI_{I_d}(E_{d,\alpha})\right) 
					\leq \left( 1+\frac{\zeta(\alpha)}{2^{\alpha-1}} \right)^d - 1
			\end{gather}
			for all $d\in\N$ and $\alpha>1$.
			Consequently, 
			\begin{gather*}
					\lim_{\alpha \nach \infty} e\!\left(N^*,d;\mathrm{Int}_d,\SI_{I_d}(E_{d,\alpha})\right)  = 0
					%, \qquad d\in\N.
			\end{gather*}
			for all $d\in\N$.
\end{theorem}
Before we give the proof of this theorem in \autoref{sect:proof} we add some comments on this result in the next section.

\section{Discussion and further results}
\begin{itemize}
\item Note that, in particular, \link{bound} yields that for every dimension~$d$ the initial error of the integration problem under consideration equals $1$. 
Thus the problem is well-scaled and we do not need to distinguish between the absolute and the normalized error criteria.
\item We stress that due to results of Smolyak and Bakhvalov the choice of linear, non-adaptive cubature rules~$\A$ in \link{quadrule} is indeed without loss of generality.
For details and further references see, e.g., \cite[Remark~1]{SW97} and Novak and Wo{\'z}niakowski~\cite[Section~4.2.2]{NW08}.
\item Moreover, observe that due to \link{upperbound} the assertion stated in \link{bound} 
can not be extended to $N\geq N^*$ provided that the smoothness $\alpha$ is sufficiently large. 
That means that at least in this case our result is sharp.
\item Note that \link{bound} can be reformulated equivalently in terms of the information complexity:
\begin{gather*}
		n\!\left(\epsilon,d;\mathrm{Int}_d,\SI_{I_d}(E_{d,\alpha})\right) \geq (\# I_d+1) \cdot 2^{d-\# I_d}
		\qquad \text{for all} \quad d\in\N \quad \text{and every} \quad \epsilon\in(0,1).
\end{gather*}
Therefore our \autoref{thm} indeed generalizes \cite[Theorem~1]{SW97} since the latter is contained as the special case $I_d=\leer$ for every $d\in\N$.
Furthermore, observe that in any case the right-hand side of the last inequality is lower bounded by $d+1$.
Hence, even for the fully permutation-invariant problem, where we have $I_d=\{1,\ldots,d\}$ for all $d\in\N$, the information complexity grows at least linearly with dimension $d$ provided that $\epsilon<1$. 
Together with some obvious estimates this proves the following corollary.
\begin{corollary}\label{cor}
		Assume $\alpha>1$ and let $b_d := d-\# I_d$.
		We study $\mathrm{Int}=(\mathrm{Int}_d)_{d\in\N}$ 
		for the sequence of $I_d$-permutation-invariant subspaces $(\SI_{I_d}(E_{d,\alpha}))_{d\in\N}$ in the worst case setting:
		\begin{itemize}
				\item If the problem is \emph{polynomially tractable} with the constants $C,p,q$ then $q \geq 1$ and
							$(b_d)_{d\in\N} \in \0(\ln d)$.
							In particular, the problem is \emph{never strongly polynomially tractable}.
				\item If the problem is \emph{uniformly weakly tractable} then $(b_d)_{d\in\N} \in o(d^{\,t})$ for all $t\in(0,1]$.
				\item If the problem is \emph{$(s,t)$-weakly tractable} for some $s,t\in(0,1]$ then 
							$(b_d)_{d\in\N} \in o(d^{\,t})$.\\
							In particular, \emph{weak tractability} implies $(b_d)_{d\in\N}\in o(d)$.
				\item If $(b_d)_{d\in\N}\notin o(d)$ then we have the \emph{curse of dimensionality}.
							In turn, already the absence of the curse implies $(b_d)_{d\in\N}\in o(d)$.
		\end{itemize}
\end{corollary}
\item For some applications it might be useful to impose
permutation-invariance conditions with respect to finitely many 
disjoint subsets $I_d^{(1)},\ldots,I_d^{(R)}\subseteq\{1,\ldots,d\}$ of the coordinates, $R>1$; see~\cite{W12b} for details.
In this case the respective subspace
\begin{gather*}
			\SI_{I_d^{(1)},\ldots,I_d^{(R)}}(E_{d,\alpha})
			:= \bigcap_{r=1}^R \SI_{I_d^{(r)}}(E_{d,\alpha}) \subset E_{d,\alpha}
\end{gather*}
consists of all $f\in E_{d,\alpha}$ which satisfy~\link{def:sym} 
for all $\sigma \in \bigcup_{r=1}^R \S_{I_d^{(r)}}$, where $\S_{I_d^{(r)}}$ 
and $\SI_{I_d^{(r)}}(E_{d,\alpha})$ are defined as before.
It turns out that \autoref{thm} remains valid for this case when we replace
the definition of $N^*$ given in~\link{def:N} by
\begin{gather*}
			N^* 
			:= N^*\!\left(d,I_d^{(1)},\ldots,I_d^{(R)}\right)
			:= \left( \prod_{r=1}^R \left(\# I_d^{(r)} + 1\right) \right) \cdot 2^{d- \sum_{r=1}^R \# I_d^{(r)}}.
\end{gather*}
Without going into details we mention that our proof given in \autoref{sect:proof} below
can be transferred almost literally to this case.
Consequently, an analogue of \autoref{cor} remains valid if we set $b_d:=d- \sum_{r=1}^R \# I_d^{(r)}$
for $d\in\N$.
\item During the last two decades many numerical problems such as (high-dimensional) integration were proven to be computationally hard.
Fortunately, it turned out that the introduction of \emph{weights} to the norm of the underlying source spaces can dramatically reduce the complexity of those problems such that they can be handled efficiently; see, e.g., \cite{NW08}.
In~\cite{W09} Wo{\'z}niakowski considered the integration problem defined above
for weighted Korobov-type spaces 
\begin{gather}\label{def_weightedE}
		E_{d,\alpha}^{\,\gamma}
		:= \left\{ f \in L_1([0,1]^d) \sep \norm{f}_\gamma := \norm{f\sep E_{d,\alpha}^{\, \gamma}} := \sup_{{\bm{k}}\in\Z^d} \abs{\widehat{f}({\bm{k}})} \frac{ \left(\overline{k_1}\cdot\ldots\cdot\overline{k_d}\right)^\alpha }{ \sqrt{ \gamma_{d,\fu(\bm{k})} } } < \infty \right\},
\end{gather}
where $d\in\N$ and $\alpha>1$, as well as $\gamma:=\{\gamma_{d,\fu(\bm{k})} \sep \bm{k}\in\Z^d\}$.
Here the quantities $\gamma_{d,\fu(\bm{k})} := \prod_{m \in \fu(\bm{k})} \gamma_{d,m}$ are \emph{product weights} with
\begin{gather*}
		\fu(\bm{k}) := \left\{ m \in\{1,\ldots,d\} \sep k_m \neq 0 \right\}
		\quad \text{and} \quad
		1 = \gamma_{d,\leer} \geq \gamma_{d,1} \geq \ldots \geq \gamma_{d,d} \geq 0
\end{gather*}
for $d\in\N$ and $\bm{k}=(k_1,\ldots,k_d)\in\Z^d$.
This generalizes \cite{SW97} since if $\gamma_{d,m}\equiv 1$ for all $m=1,\ldots,d$ and $d\in\N$ then
$E_{d,\alpha}^{\,\gamma}=E_{d,\alpha}$.
Theorem~1 in~\cite{W09} states that
\begin{gather*}
		\eta_{d,N} \leq e(N,d; \mathrm{Int}_d, E_{d,\alpha}^{\,\gamma}) \leq 1 
		\quad \text{for all} \quad \alpha>1, \quad d\in\N, \quad \text{and} \quad N<2^d,
\end{gather*}
where $\eta_{d,N}$ denotes the $(N+1)$st largest weight in the sequence $(\gamma_{d,\fu})_{\leer\subseteq\fu\subseteq\{1,\ldots,d\}}\subset [0,1]$.
So what about weights in our setting?

First observe that the subset of $I_d$-permutation-invariant functions in $E_{d,\alpha}^{\,\gamma}$ again forms a linear subspace. In accordance with our previous notation we denote it by $\SI_{I_d}(E_{d,\alpha}^{\,\gamma})$.
Furthermore, we note that the Fourier coefficients of an $I_d$-permutation-invariant function $f$ are $I_d$-permutation-invariant again, i.e., they satisfy $\widehat{f}(\bm{k}) = \widehat{f}(\sigma(\bm{k}))$ for every $\bm{k}\in\Z^d$ and all $\sigma\in\S_{I_d}$.
Since clearly also $\left( \overline{k_1}\cdot\ldots\cdot\overline{k_d}\right)^\alpha$ is $I_d$-permutation-invariant, $f\in E_{d,\alpha}^{\,\gamma}$ belongs to the unit ball of the subspace $\SI_{I_d}(E_{d,\alpha}^{\,\gamma})$ if and only if
\begin{gather*}
		\abs{\widehat{f}(\bm{k})} 
		\leq \min_{\sigma\in\S_{I_d}} \frac{\sqrt{\gamma_{d,\fu(\sigma(\bm{k}))}}}{\left( \overline{k_1}\cdot\ldots\cdot\overline{k_d}\right)^\alpha}
		=: \frac{\sqrt{\mu_{d,\fu(\bm{k})}}}{\left( \overline{k_1}\cdot\ldots\cdot\overline{k_d}\right)^\alpha},
		\quad \bm{k}\in\Z^d,
\end{gather*}
where the new weights $\mu_{d,\fu}$ are defined as a product of an order-dependent part (\wrt coordinates from $I_d$) and a part consisting of usual product weights. For the ease of notation, in what follows we assume permutation-invariance \wrt the first $\# I_d$ coordinates, i.e., $I_d=\{1,2,\ldots,\#I_d\}$. Then
\begin{align}\label{weights_mu}
		\mu_{d,\fu(\bm{k})}
		= \min_{\sigma\in\S_{I_d}} \prod_{m\in \fu_{I_d}(\sigma(\bm{k}))} \gamma_{d,m} \cdot \prod_{m\in \fu_{I_d^c}(\bm{k})} \gamma_{d,m}
		= \prod_{i=1}^{\#\fu_{I_d}(\bm{k})} \gamma_{d,\#I_d-i+1} \cdot \prod_{m\in \fu_{I_d^c}(\bm{k})} \gamma_{d,m}.
\end{align}
Here we set $\fu_{I_d}(\bm{k}):=\fu(\bm{k}) \cap I_d$ and $\fu_{I_d^c}(\bm{k}):=\fu(\bm{k}) \setminus I_d$ to denote the
support of $\bm{k}\in\Z^d$ \wrt the sets $I_d$ and $I_d^c = \{1,\ldots,d\}\setminus I_d$, respectively.
Thus we have 
\begin{equation*}
	\SI_{I_d}(E_{d,\alpha}^{\,\gamma})=\SI_{I_d}(E_{d,\alpha}^{\, \mu}),
\end{equation*}
where $E_{d,\alpha}^\mu$ is defined by \link{def_weightedE} with $\gamma_{d,\fu(\bm{k})}$ replaced by $\mu_{d,\fu(\bm{k})}$ and $\SI_{I_d}$ once more denotes the restriction to $I_d$-permutation-invariant functions.

Again the lower bound for the error of numerical integration depends on the $(N+1)$th largest weight~$\mu_{d,\fu(\bm{k})}$ that can appear.
To formalize this point let
\begin{gather*}
		\psi \colon \{0,1,\ldots, \# \nabla_d - 1 \} \nach \nabla_d
\end{gather*}
denote a rearrangement of the multi-indices $\bm{k}$ from the set
\begin{gather}\label{nabla}
		\nabla_d 
		:= \left\{ {\bm{k}}=(\bm{m},\bm{l}) \in \{0,1\}^{\#I_d} \times \{0,1\}^{(d-\# I_d)} \sep m_1\leq \ldots \leq m_{\# I_d}\right\} \subset \{0,1\}^d \subset \Z^d
\end{gather}
such that the corresponding weights $\mu_{d,\fu(\bm{k})}$ possess a non-increasing ordering, \ie
\begin{gather}\label{ordering}
		\nu_{d,n} := \mu_{d,\fu(\bm{\psi(n))}} \geq \mu_{d,\fu(\bm{\psi(n+1)})} = \nu_{d,n+1} \quad \text{for all} \quad n=0,1,\ldots,\#\nabla_d-2.
\end{gather}
Then the weighted analogue of \autoref{thm} reads as follows:
\begin{theorem}\label{thm_weighted}
			Assume $\alpha>1$, for every $d\in\N$ let $I_d \subseteq \{1,\ldots,d\}$ be given, and consider $N^*(d,I_d)$ defined as in \link{def:N}.
			Then we have
			\begin{gather}\label{bound_weighted}
\nu_{d,N} \leq	e\!\left(N,d;\mathrm{Int}_d,\SI_{I_d}(E_{d,\alpha}^{\,\mu})\right) \leq 1
			\end{gather}
			for all $d\in\N$ and every $N < N^*(d,I_d)$.
\end{theorem}

We illustrate \autoref{thm_weighted} by two examples.
\begin{example}
	Consider the fully permutation-invariant problem, i.e., assume $I_d=\{1,\ldots,d\}$ for every $d\in\N$. For simplicity let
	\begin{gather*}
		\gamma_{d,m} := \frac{1}{m} \quad \text{for every} \quad m\in\{1,\ldots,d\} \quad \text{and all} \quad d\in\N.
	\end{gather*}
	Then $\nabla_d=\{\bm{k}\in\{0,1\}^d \sep k_1\leq\ldots\leq k_d\}$ and $\# \nabla_d=d+1$. Moreover, we conclude $\mu_{d,\fu(\bm{k})}=\prod_{i=1}^{\# \fu(\bm{k})} 1/(d-i+1)$, $\bm{k}\in\nabla_d$, and  consequently
	\begin{align*}
		\nu_{d,0} &= 1,\\
		\nu_{d,1} &= 1/d, \\
		\nu_{d,2} &= 1/(d\cdot (d-1)), \\
		&\;\vdots\\
		\nu_{d,d} &= 1/(d!).
	\end{align*}
	Hence, for $\alpha>1$ and $d\in\N$
	we obtain $e(0,d;\mathrm{Int}_d; \SI_{I_d}(E^{\,\mu}_{d,\alpha})) = 1$ and
	\begin{gather*}
		e(N,d;\mathrm{Int}_d; \SI_{I_d}(E^{\,\mu}_{d,\alpha})) \geq \frac{1}{d \cdot (d-1) \cdot \ldots \cdot (d-N+1)},
	\end{gather*}
	if $N \in\{ 1,\ldots,d \}$.\hfill $\square$
\end{example}

Our second, more sophisticated example shows that permutation-invariance is not as powerful as additional knowledge modeled by weights. It generalizes an assertion due to Wo\'{z}niakowski~\cite[p. 648]{W09}.
\begin{example}
	For $d\in\N$ let $\#I_d \in \{0,1,\ldots,d\}$ be arbitrary and assume permutation-invariance \wrt the first $\#I_d$ coordinates in dimension $d$. Furthermore, let
	\begin{gather*}
		\gamma_{d,1}:=\ldots:=\gamma_{d,\#I_d}:=1 \quad \text{and} \quad 1 \geq \gamma_{d,\#I_d+1} \geq \ldots \geq \gamma_{d,d} > 0.
	\end{gather*}
	That is, we assume the weights only act on coordinates without permutation-invariance conditions and leave the remaining coordinates unweighted.
Thus we ask how much permutation-invariance can relax (by now well-established) necessary conditions on the sequences $(\gamma_{d,m})_{m=1}^d$, $d\in\N$, for (strong) polynomial tractability:

	Assuming polynomial tractability with constants $C\geq 1$, $p>0$ and $q\geq 0$ it is easy to check that
	\begin{gather*}
		e(N,d;\mathrm{Int}_d, \SI_{I_d}(E^{\,\mu}_{d,\alpha})) \leq (2\, C)^{1/p} \, d^{q/p} \, N^{-1/p},
	\end{gather*}
	see \cite{W09} for details.
	For every $d\in\N$ and all $\kappa>p$ this implies
	\begin{gather*}
		\sum_{N=0}^{\infty} e(N,d;\mathrm{Int}_d, \SI_{I_d}(E^{\,\mu}_{d,\alpha}))^{\kappa} \leq 1 + (2C)^{\kappa/p}\, d^{q\, \kappa/p} \, \zeta(\kappa/p) \leq C' \, d^{q\, \kappa/p}
	\end{gather*}
	with some finite constant $C'>0$ which only depends on $C,\kappa$, and $p$.
	On the other hand, \autoref{thm_weighted} provides the lower bound
	\begin{align*}
		\sum_{N=0}^{\infty} e(N,d;\mathrm{Int}_d; \SI_{I_d}(E^{\,\mu}_{d,\alpha}))^{\kappa}
		&\geq \sum_{N=0}^{N^*(d,I_d)-1} \nu_{d,N}^{\kappa} 
		= \sum_{N=0}^{\#\nabla_d-1} \mu_{d,\fu(\bm{\psi(N)})}^{\kappa}
		= \sum_{\bm{k}\in\nabla_{d}} \mu_{d,\fu(\bm{k})}^{\kappa} \\
		\qquad &= \sum_{\substack{\bm{m}\in\{0,1\}^{\#I_d} \\ m_1\leq\ldots \leq m_{\#I_d}}} 1 \cdot \sum_{\bm{l}\in\{0,1\}^{d-\#I_d}} \left( \prod_{j\in\fu(\bm{l})} \gamma_{d,\#I_d+j} \right)^{\kappa} \\
		&= (\#I_d+1) \cdot \prod_{j=1}^{d-\#I_d} \left( 1 + \gamma_{d,\#I_d+j}^\kappa \right).
	\end{align*}
	In conclusion
	\begin{gather}\label{cond}
		\limsup_{d\nach\infty} \frac{(\#I_d+1) \cdot \prod_{j=1}^{d-\#I_d} \left( 1 + \gamma_{d,\#I_d+j}^\kappa \right)}{d^{q\,\kappa/p}} < \infty
	\end{gather}
	is a necessary condition for (strong) polynomial tractability.
	In particular, this implies that we can have strong polynomial tractability ($q=0$) only if the number of permutation-invariant coordinates $\#I_d$ is uniformly bounded. If so, then \link{cond} yields that
	\begin{gather*}
		\limsup_{d\nach\infty} \sum_{j=1}^{d} \gamma_{d,j}^\kappa < \infty 
		\quad \text{for all} \quad 
		\kappa > p
	\end{gather*}
	which resembles the assertion stated in~\cite[Cor.~1]{W09}.
	If we assume polynomial tractability with $q>0$ then we find that for all $\kappa>p$
	\begin{gather*}
		\#I_d \in \0(d^{q \, \kappa/p}) 
		\quad \text{and} \quad
		\limsup_{d\nach\infty} \frac{\sum_{j=1}^{d-\#I_d} \gamma_{d,\#I_d+j}^\kappa}{\ln d} < \infty.
	\end{gather*}
	Note that if $\# I_d$ is bounded then the latter condition again coincides with the known condition given by 
	Wo\'{z}niakowski. Otherwise, if $\#I_d$ grows like say $d^\beta$ for some $\beta\in(0,1]$, then our condition turns out to be less restrictive, i.e., we may have polynomial tractability although polynomially many coordinates are unweighted (but permutation-invariant).\hfill $\square$
\end{example}
\end{itemize} 

\section{Proofs}\label{sect:proof}
In order to prove \autoref{thm} we basically combine the ideas stated in \cite{SW97,W09} with
the technique developed in \cite{W12b}. 
Since the proof is a little bit technical we divide it into three steps which are
organized as follows.

In a first step we show that for any given integration rule that uses $N<N^*$ function values 
there exists a certain fooling function which shows that 
$e(N,d;\mathrm{Int}_d,\SI_{I_d}(E_{d,\alpha})) \geq 1$.
Afterwards, in Step 2, we notice that this lower bound is sharp, because
$e(0,d;\mathrm{Int}_d,\SI_{I_d}(E_{d,\alpha}))\leq 1$ for every $d\in\N$ and all $\alpha>1$.
Finally, we present a cubature rule that uses at most $N=N^*$ function values, 
whereas its worst case error is no larger than the bound stated in \link{upperbound}.

\begin{proof}[Proof (\autoref{thm})]
\emph{Step 1.}
Following Sloan and Wo{\'z}niakowski~\cite{SW97} we fix $\alpha>1$, as well as $d\in\N$, and consider
an arbitrary linear cubature rule $\A:=\A_{N,d}$, given by \link{quadrule},
with $N<N^*(d,I_d):=(\# I_d + 1)\cdot 2^{d-\# I_d}$.
Without loss of generality let us assume that $I_d=\{1,\ldots,\#I_d\}$.
We will show that there exists a function $f_N$ in the unit ball of $\SI_{I_d}(E_{d,\alpha})$ 
such that $\A(f_N)=0$, whereas the integral of $f_N$ equals $1$.

For ${{\bm{k}}}\in\Z^d$ let $e_{{\bm{k}}} := \exp(2\pi \i {{\bm{k}}}\, \cdot)$. 
Following the lines of \cite{W12b} we define a linear operator $\SI_{I_d}\colon E_{d,\alpha} \nach E_{d,\alpha}$ called \emph{symmetrizer} by
\begin{gather*}
		(\SI_{I_d}e_{\bm{k}})(\bm{x}) := \frac{1}{\#\S_{I_d}} \sum_{\sigma \in \S_{I_d}} e_{{\bm{k}}}(\sigma(\bm{x})), \qquad \bm{x}\in[0,1]^d,
\end{gather*}
and continuous extension.
Therein $\S_{I_d}$ is given in \link{def:permut_set}.
Then $\SI_{I_d} e_{\bm{k}} = (1/\#\S_{I_d}) \sum_{\lambda\in\S_{I_d}} e_{\lambda^{-1}({\bm{k}})}(\cdot)$ is $I_d$-permutation-invariant in the sense of \link{def:sym}, i.e.,
$(\SI_{I_d} e_{\bm{k}})(\sigma(\bm{x})) = (\SI_{I_d} e_{\bm{k}})(\bm{x})$
for all $\bm{x}\in[0,1]^d$ and any $\sigma\in\S_{I_d}$.
Moreover, for ${\bm{k}}\in\Z^d$ let $M_{I_d}({\bm{k}})$ be defined as in \cite{W12b}. 
That is, $M_{I_d}({\bm{k}})!$ denotes the number of different permutations $\sigma\in\S_{I_d}$ such that $\sigma({\bm{k}})=\lambda(\bm{k})$ for any fixed $\lambda\in\S_{I_d}$.

To prove the claim we consider the set $\nabla_d$ introduced in \link{nabla}.
Observing that $\# \nabla_d=N^*(d,I_d)=(\# I_d + 1)\cdot 2^{d-\# I_d}$
we choose a bijection
\begin{gather}\label{psi_arb}
			\psi \colon \{0,1,\ldots, N^*(d,I_d)-1\} \nach \nabla_d.
\end{gather}

Furthermore, we consider the homogeneous linear system
\begin{gather*}
		\sum_{n=0}^{N} a_{n} \cdot \frac{(\SI_{I_d} e_{\bm{\psi(n)}})(\bm{t}^{(i)})}{M_{I_d}(\bm{\psi(n)})!} 
		= 0, 		
		\qquad i=1,\ldots,N,
\end{gather*}
which consists of $N<N^*(d,I_d)$ linear equations in $N+1$ complex variables $a_n$, $n=0,1,\ldots,N$.
Here the points $\bm{t}^{(i)}\in[0,1]^d$, $i=1,\ldots,N$, denote the integration nodes used by the cubature rule~$\A$ applied to the function $0\in\SI_{I_d}(E_{d,\alpha})$.
Clearly, we can select a non-trivial solution $\bm{a}=(a_n)_{n=0}^N\in\C^{N+1}$ of this system, scaled such that 
for some $n^*\in\{0,1,\ldots,N\}$
\begin{gather*}
			a_{n^*} = 1 \geq \max_{n=0,1,\ldots,N} \abs{a_{n}}.
\end{gather*}
Next we define the function $f_N \colon [0,1]^d\nach\C$ by
\begin{gather}\label{def:f_N}
			f_N(\bm{x}) 
			:= \#\S_{I_d} \cdot (\SI_{I_d} e_{-\bm{\psi(n^*)}})(\bm{x}) \cdot \left( \sum_{n=0}^N a_{n} \cdot \frac{(\SI_{I_d} e_{\bm{\psi(n)}})(\bm{x})}{M_{I_d}(\bm{\psi(n)})!} \right), \qquad \bm{x}\in[0,1]^d.
\end{gather}
Observe that then $f_N(\bm{t}^{(i)})=0$ for all $i=1,\ldots,N$ and thus we have $\A(f_N)=0$.
Moreover, we see that $f_N$ is $I_d$-permutation-invariant as a product of $I_d$-permutation-invariant functions.
It remains to show that $f_N$ is a suitable fooling function for our integration problem.
That is, we show that $f_N$ is an element of the unit ball of $\SI_{I_d}(E_{d,\alpha})$, 
i.e., its Fourier coefficients satisfy
\begin{gather}\label{unitball}
			\abs{\widehat{f_N}({\bm{k}})} \leq \frac{1}{\left( \overline{k_1}\cdot\ldots\cdot\overline{k_d}\right )^\alpha} \quad \text{for every} \quad \bm{k}=(k_1,\ldots,k_d)\in\Z^d,
\end{gather} 
and that its integral is as large as possible, i.e., $\mathrm{Int}_d(f_N)=\widehat{f_N}(\bm{0})=1$.

For ${\bm{k}}\in\Z^d$ we calculate
\begin{align*}
		\widehat{f_N}(\bm{k}) 
		&= \distr{f_N}{e_{\bm{k}}}_{L_2} \\
		&= \#\S_{I_d} \distr{\left(\frac{1}{\#\S_{I_d}} \sum_{\sigma\in\S_{I_d}} e_{\sigma^{-1}(-\bm{\psi(n^*)})}\right) \left( \sum_{n=0}^N \frac{a_{n}}{M_{I_d}(\bm{\psi(n)})!} \, \frac{1}{\#\S_{I_d}} \sum_{\lambda \in \S_{I_d}} e_{\lambda^{-1}(\bm{\psi(n)})} \right)}{e_{\bm{k}}}_{L_2} \\
		&= \frac{1}{\#\S_{I_d}} \sum_{\sigma\in\S_{I_d}} \sum_{\lambda\in\S_{I_d}} \sum_{n=0}^N \frac{a_n}{M_{I_d}(\bm{\psi(n)})!} \distr{e_{-\sigma(\bm{\psi(n^*)})} e_{\lambda(\bm{\psi(n)})}}{e_{\bm{k}}}_{L_2}
\end{align*}
and
\begin{align*}
		\distr{e_{-\sigma(\bm{\psi(n^*)})} e_{\lambda(\bm{\psi(n)})}}{e_{\bm{k}}}_{L_2} 
		&= \int_{[0,1]^d} \exp\!\bigg( 2\pi \i \, \big(\lambda(\bm{\psi(n)})-{\bm{k}}-\sigma(\bm{\psi(n^*)})\big) \bm{x} \bigg) \d\bm{x} \\
		&= \delta_{\left[\lambda(\bm{\psi(n)})-{\bm{k}}-\sigma(\bm{\psi(n^*)})=0\right]},
\end{align*}
where $\delta_{[C]}$ equals $1$ if the condition $C$ is fulfilled and $0$ otherwise.
In particular from $\widehat{f_N}(\bm{k})\neq 0$ it follows that there exist $\lambda$ and $\sigma$ in $\S_{I_d}$, as well as $n\in\{0,1,\ldots,N\}$, such that
\begin{gather}\label{relation}
		\bm{k} = \lambda(\bm{\psi(n)})-\sigma(\bm{\psi(n^*)} \in \{-1,0,1\}^d,
\end{gather}
since $\bm{\psi(n)}\in\nabla_d\subset \{0,1\}^d$ for all $n$ and for every permutation $\lambda\in \S_{I_d}$ the multi-index $\lambda(\bm{\psi(n)})$ is an element of $\{0,1\}^d$, too.

Hence we arrive at $\widehat{f_N}(\bm{k})=0$ if $\bm{k}\in\Z^d\setminus\{-1,0,1\}^d$ and
\begin{gather*}
		\widehat{f_N}({\bm{k}}) =
		\frac{1}{\#\S_{I_d}} \sum_{\sigma\in\S_{I_d}} \delta_{\left[ {\bm{k}}+\sigma(\bm{\psi(n^*)}) \in \{0,1\}^d \right]} \cdot \sum_{\lambda\in\S_{I_d}} \left( \sum_{n=0}^N \frac{a_{n}}{M_{I_d}(\bm{\psi(n)})!} \cdot \delta_{\left[\lambda(\bm{\psi(n)})={\bm{k}}+\sigma(\bm{\psi(n^*)})\right]} \right)
\end{gather*}
if $\bm{k}\in\{-1,0,1\}^d$.

We will show that the summation within the brackets can be reduced to 
at most one non-vanishing term.
Therefore assume $\sigma\in\S_{I_d}$ to be fixed and consider
\begin{gather*}
		\bm{h}:=\bm{k}+\sigma(\bm{\psi(n^*)})\in\{0,1\}^d.
\end{gather*}
For this $\bm{h}$ there exists one (and only one) multi-index $\bm{j}\in\nabla_d$ such that
\begin{gather*}
		(j_m)_{m\in I_d^c} = (h_m)_{m\in I_d^c} \quad \text{and} \quad \#\{m\in I_d \sep j_m=1\} = \#\{m\in I_d \sep h_m=1\},
\end{gather*}
because $\lambda$ and $\sigma$ in $\S_{I_d}$ leave the coordinates in $I_d^c=\{\#I_d+1,\ldots,d\}$ fixed.
Since the mapping $\psi$ was assumed to be bijective, $\bm{j}$ uniquely defines
\begin{gather*}
		n(\bm{k},n^*,\sigma) := \psi^{-1}(\bm{j}) \in \{0,1,\ldots,\#\nabla_d-1\}.
\end{gather*}
Hence, there can be at most one $n\in\{0,1,\ldots,N\}$ with $\lambda(\bm{\psi(n)})=\bm{h}$.
If so, then there exist exactly $M_{I_d}(\bm{\psi(n(\bm{k},n^*,\sigma))})!$ different permutations $\lambda\in\S_{I_d}$ such that
\begin{gather*}
		\lambda(\bm{\psi(n(\bm{k},n^*,\sigma))})=\bm{h}=\bm{k}+\sigma(\bm{\psi(n^*)}).
\end{gather*}
If $\bm{k}=\bm{0}\in\Z^d$ then $\bm{h}=\sigma(\bm{\psi(n^*)})$ which implies $\bm{j}=\bm{\psi(n^*)}$ and thus
$n(\bm{0},n^*,\sigma)=\psi^{-1}(\bm{j})=n^*$ for all $\sigma\in\S_{I_d}$.

Consequently, we obtain
\begin{gather}\label{fourier}
		\widehat{f_N}(\bm{k}) =
		\frac{1}{\#\S_{I_d}} \sum_{\sigma\in\S_{I_d}} \delta_{\left[ \bm{k}+\sigma(\bm{\psi(n^*)}) \in \{0,1\}^d \right]} \cdot 
		\delta_{\left[ n(\bm{k},n^*,\sigma)\leq N \right] } \cdot a_{n(\bm{k},n^*,\sigma)},
		\quad \bm{k}\in\{-1,0,1\}^d,
\end{gather}
which particularly yields
\begin{gather*}
		\mathrm{Int}_d(f_N) 
		= \widehat{f_N}(\bm{0}) 
		= \frac{1}{\#\S_{I_d}} \sum_{\sigma\in\S_{I_d}} a_{n^*} 
		= a_{n^*} 
		= 1.
\end{gather*}
Furthermore, formula \link{fourier} implies
\begin{align*}
		\abs{\widehat{f_N}(\bm{k})} \, \left( \overline{k_1}\cdot\ldots\cdot\overline{k_d}\right )^\alpha
		&= \abs{\widehat{f_N}(\bm{k})} \\
		&\leq \frac{1}{\#\S_{I_d}} \sum_{\sigma\in\S_{I_d}} \abs{ \delta_{\left[ \bm{k}+\sigma(\bm{\psi(n^*)}) \in \{0,1\}^d \right]} }\cdot 
		\abs{ \delta_{\left[ n(\bm{k},n^*,\sigma)\leq N \right] }} \cdot \abs{a_{n(\bm{k},n^*,\sigma)}} \\
		&\leq \frac{1}{\#\S_{I_d}} \sum_{\sigma\in\S_{I_d}} a_{n^*} = 1,
\end{align*}
for $\bm{k}\in\{-1,0,1\}^d$.
Together with
$\abs{\widehat{f_N}(\bm{k})} \, \left( \overline{k_1}\cdot\ldots\cdot\overline{k_d}\right )^\alpha=0$ for $\bm{k}\in\Z^d \setminus \{-1,0,1\}^d$ this finally proves \link{unitball} and completes this step.

\emph{Step 2}.
Clearly, we have $e(0,d;\mathrm{Int}_d,\SI_{I_d}(E_{d,\alpha}))\leq 1$ 
for all $d\in\N$ and any $I_d\subseteq\{1,\ldots,d\}$
since the worst case error of the zero algorithm $\A_{0,d}\equiv 0$ is given by
\begin{gather*}
		\Delta^{\mathrm{wor}}\! \left(\A_{0,d}; \mathrm{Int}_d, \SI_{I_d}(E_{d,\alpha})\right) 
		= \sup_{\norm{f \sep \SI_{I_d}(E_{d,\alpha})} \leq 1} \abs{\widehat{f}(\bm{0})} 
		\leq 1.
\end{gather*}

\emph{Step 3}.
To show \link{upperbound} we once more follow the arguments given in~\cite{SW97}.
There it has been shown that the $2^d$-point product-rectangle rule
\begin{gather*}
	 	\algR_{2^d,d} \colon E_{d,\alpha} \nach \C, \qquad
		f \mapsto \algR_{2^d,d}(f) := \frac{1}{2^d} \sum_{{\bm{j}} \in \{0,1\}^d} f \!\left( \frac{j_1}{2}, \frac{j_2}{2},\ldots,\frac{j_d}{2} \right),
\end{gather*}
is a suitable algorithm for $\mathrm{Int}_d$ on $E_{d,\alpha}$.
For $d\in\N$ and $\alpha>1$ its worst case error was found to be 
\begin{gather}\label{wce}
		\Delta^{\mathrm{wor}}(\algR_{2^d,d}; \mathrm{Int}_d, E_{d,\alpha}) 
		= \left( 1+\frac{\zeta(\alpha)}{2^{\alpha-1}} \right)^d - 1.
\end{gather}
Since $\SI_{I_d}(E_{d,\alpha})$ is a linear subspace of 
$E_{d,\alpha}$, where $\norm{\cdot \sep \SI_{I_d}(E_{d,\alpha})} = \norm{\cdot \sep E_{d,\alpha}}$, 
we can restrict the algorithm $\algR_{2^d,d}$ to $I_d$-permutation-invariant functions in $E_{d,\alpha}$
and its worst case error remains bounded from above by~\link{wce}.
On the other hand, it is easy to see that
\begin{align}
		\A_{N^*,d}(f) 
		&:= \sum_{{\bm{k}}\in\nabla_d} \frac{\# \S_{I_d}}{2^d \cdot M_{I_d}({\bm{k}})!} \cdot f \!\left( \frac{{\bm{k}}}{2} \right) \label{eq:algoA}\\
		&= \frac{1}{2^d} \sum_{{\bm{k}}\in\nabla_d} \sum_{\sigma \in \S_{I_d}} \frac{f(\sigma({\bm{k}})/2)}{M_{I_d}({\bm{k}})!} 
		= \frac{1}{2^d} \sum_{{\bm{j}} \in \{0,1\}^d} f \!\left( \frac{{\bm{j}}}{2} \right)
		= \algR_{2^d,d}(f) \nonumber
\end{align}
on $\SI_{I_d}(E_{d,\alpha})$.
Hence the restriction of $\algR_{2^d,d}$, i.e., the algorithm $\A_{N^*,d}$ defined by~\link{eq:algoA},
is a cubature rule in the sense of~\link{quadrule} that uses not more than $\# \nabla_d = N^*=N^*(d,I_d)$ 
function values in dimension $d$; see~\link{nabla}.
Consequently,
\begin{align*}
		e(N^*,d;\mathrm{Int}_d,\SI_{I_d}(E_{d,\alpha})) 
		&\leq \Delta^{\mathrm{wor}}(\A_{N^*,d}; \mathrm{Int}_d, \SI_{I_d}(E_{d,\alpha})) \\
		&\leq \Delta^{\mathrm{wor}}(\algR_{2^d,d}; \mathrm{Int}_d, E_{d,\alpha})
		= \left( 1+\frac{\zeta(\alpha)}{2^{\alpha-1}} \right)^d - 1
\end{align*}
which completes the proof.
\end{proof}

The proof of \autoref{thm_weighted} can be derived using only one additional argument.
As in the previous proof, we construct a suitable fooling function $g_N$ that lies in the unit ball of $\SI_{I_d}(E_{d,\alpha}^{\,\mu})$ and for which $\mathrm{Int}_d(g_N)$ is as large as possible while $\A(g_N)=0$.
\begin{proof}[Proof (\autoref{thm_weighted})]
Since we now deal with the weighted case we need to show that for every $\bm{k}=(k_1,\ldots,k_d)\in\Z^d$ we have
\begin{gather*}
			\abs{\widehat{g_N}({\bm{k}})} 
			\leq \frac{\sqrt{\mu_{d,\fu(\bm{k})}}}{\left( \overline{k_1}\cdot\ldots\cdot\overline{k_d}\right )^\alpha}
\end{gather*}
instead of \link{unitball}.
We start by selecting a bijection $\psi\colon \{0,1,\ldots,N^*(d,I_d)-1\}\nach\nabla_d$ which provides a non-increasing ordering of the weights; see \link{psi_arb} and \link{ordering}.
Setting
\begin{gather*}
			g_N(\bm{x}):=\sqrt{\mu_{d,\fu(\bm{\psi(N)})} \cdot \mu_{d,\fu(\bm{\psi(n^*)})}} \cdot f_N(\bm{x}), \qquad \bm{x}\in[0,1]^d,
\end{gather*}
it now remains to show that for every $\bm{k}\in\Z^d$ with $\widehat{g_N}(\bm{k})\neq 0$
we can estimate
\begin{gather}\label{weight_est}
			\mu_{d,\fu(\bm{\psi(N)})} \cdot \mu_{d,\fu(\bm{\psi(n^*)})} 
			\leq \mu_{d,\fu(\bm{k})}.
\end{gather}
Recall that here the relation of $\bm{k}$ and $n^*$ is given by $\bm{k}=\lambda(\bm{\psi(n)})-\sigma(\bm{\psi(n^*)})$ for some $\lambda,\sigma\in\S_{I_d}$ and a certain $n\in\{0,1,\ldots,N\}$; see \link{relation}.

To prove \link{weight_est} we first note that the latter representation of $\bm{k}$ yields
\begin{gather*}
			\# \fu_{I_d}(\bm{k}) \leq \minn{\# \fu_{I_d}(\bm{\psi(n)})+ \# \fu_{I_d}(\bm{\psi(n^*)}), \# I_d},
\end{gather*}
since $\lambda$ and $\sigma$ do not change the size of the support of the respective multi-indices $\bm{\psi(n)}$ and $\bm{\psi(n^*)}$ in $\nabla_d$.
Consequently we obtain
\begin{align}\label{est_order}
			\prod_{i=1}^{\#\fu_{I_d}(\bm{k})} \gamma_{d,\#I_d-i+1} 
			&\geq \prod_{i=1}^{\#\fu_{I_d}(\bm{\psi(n)})} \gamma_{d,\#I-i+1} \cdot \prod_{i=\#\fu_{I_d}(\bm{\psi(n)})+1}^{\minn{\# \fu_{I_d}(\bm{\psi(n)})+ \# \fu_{I_d}(\bm{\psi(n^*)}), \# I_d}} \gamma_{d,\#I_d-i+1} \nonumber\\
			&= \prod_{i=1}^{\#\fu_{I_d}(\bm{\psi(n)})} \gamma_{d,\#I_d-i+1} \cdot \prod_{i=1}^{\minn{\# \fu_{I_d}(\bm{\psi(n^*)}), \# I_d-\#\fu_{I_d}(\bm{\psi(n)})}} \gamma_{d,\#I_d-\#\fu_{I_d}(\bm{\psi(n)})-i+1} \nonumber\\
			&\geq \prod_{i=1}^{\#\fu_{I_d}(\bm{\psi(n)})} \gamma_{d,\#I_d-i+1} \cdot \prod_{i=1}^{\#\fu_{I_d}(\bm{\psi(n^*)})} \gamma_{d,\#I_d-i+1},
\end{align}
by exploiting the general assumption $1\geq \gamma_{d,1}\geq \ldots \geq \gamma_{d,\#I_d} \geq 0$.

Moreover, for the coordinates related to the product weight part we conclude
\begin{gather*}
			\fu_{I_d^c}(\bm{k}) \subseteq \fu_{I_d^c}(\bm{\psi(n)}) \cup \fu_{I_d^c}(\bm{\psi(n^*)})
\end{gather*}
which immediately implies that
\begin{gather}\label{est_prod}
			\prod_{m\in \fu_{I_d^c}(\bm{k})} \gamma_{d,m}
			\geq \prod_{m\in \fu_{I_d^c}(\bm{\psi(n)})} \gamma_{d,m} \cdot \prod_{m\in \fu_{I_d^c}(\bm{\psi(n^*)})} \gamma_{d,m}.
\end{gather}
Together with the representation of $\mu$ stated in \link{weights_mu} the estimates \link{est_order} and \link{est_prod}
imply 
\begin{gather*}
			\mu_{d,\fu(\bm{\psi(n)})} \cdot \mu_{d,\fu(\bm{\psi(n^*)})} \leq \mu_{d,\fu(\bm{k})}. 
\end{gather*}
This leads us to the claimed bound \link{weight_est} using the imposed ordering $\mu_{d,\fu(\bm{\psi(N)})} \leq \mu_{d,\fu(\bm{\psi(n)})}$ for $n=0,1,\ldots,N$.

Finally, since $\widehat{f_N}(\bm{0})=1$, the integral of $g_N$ is lower bounded by
\begin{gather*}
			\mathrm{Int}_d(g_N) 
			=\widehat{g_N}(\bm{0}) 
			=\sqrt{\mu_{d,\fu(\bm{\psi(N)})} \cdot \mu_{d,\fu(\bm{\psi(n^*)})}} \cdot \widehat{f_N}(\bm{0})
			\geq \mu_{d,\fu(\bm{\psi(N)})} = \nu_{d,N}.
\end{gather*}
The fact that $\A(g_N)=\sqrt{\mu_{d,\fu(\bm{\psi(N)})} \cdot \mu_{d,\fu(\bm{\psi(n^*)})}} \cdot \A(f_N)=0$ completes the proof.
\end{proof}

\section*{Acknowledgements}
\addcontentsline{toc}{section}{Acknowledgments}
I want to thank Henryk Wo{\'z}niakowski for pointing out the existence of the reference~\cite{W09} to me.
Moreover, I like to thank Dirk Nuyens and Gowri Suryanarayana for motivating me to study integration problems of permutation-invariant functions, as well as for their kind hospitality during my stays in Leuven in October '12 and in March '13.


\addcontentsline{toc}{chapter}{References}
\begin{thebibliography}{Wo{\'z}09}

\bibitem[NW08]{NW08}
\textsc{E. Novak and H. Wo{\'z}niakowski} - \textit{Tractability of
  {M}ultivariate {P}roblems. {V}ol. {I}: {L}inear {I}nformation}. EMS Tracts in
  Mathematics~6. European Mathematical Society (EMS), Z\"urich. 2008.

\bibitem[Sie13]{P13}
\textsc{P. Siedlecki} - \textit{Uniform weak tractability}. 
  J.~Complexity \textbf{29}, \textit{2013}, pp.~438--453.
  
\bibitem[SW97]{SW97}
\textsc{I. Sloan and H. Wo{\'z}niakowski} - \textit{An intractability result
  for multiple integration}. Mathematics of Computation \textbf{66},
  \textit{1997}, pp.~1119--1124.

\bibitem[Wei12]{W12b}
\textsc{M. Weimar} - \textit{The complexity of linear tensor product problems
  in (anti)symmetric {H}ilbert spaces}. J.~Approx. Theory \textbf{164(10)},
  \textit{2012}, pp.~1345--1368.

\bibitem[Wo{\'z}09]{W09}
\textsc{H. Wo{\'z}niakowski} - \textit{Tractability of multivariate integration
  for weighted {K}orobov spaces: my 15 year partnership with {I}an {S}loan},
  in: P. L'Ecuyer and A.B. Owen (Eds.) - Monte Carlo and Quasi-Monte Carlo
  Methods 2008. Springer, Berlin. 2009, pp.~637--653.
\end{thebibliography}
\end{document}